\documentclass[12pt]{article}               
\usepackage{geometry}
\usepackage{cite}
\usepackage{multirow}
\usepackage{amsmath}
\usepackage{amssymb}
\usepackage[usenames]{color}
\usepackage{graphicx}          
 \usepackage{amsthm}                              
\usepackage{epsfig}    
 \usepackage{helvet}	
 
\newtheorem{definition}{Definition}
\newtheorem{theorem}{Theorem}

\newtheorem{remark}{Remark}
\newtheorem{problem}{Problem}

\newtheorem{assumption}{Assumption}

\title{\large Partial stabilization of nonholonomic systems\\ with application to multi-agent coordination}

\author{Victoria Grushkovskaya\thanks{
Institute of Mathematics, Alpen-Adria University of Klagenfurt,
       9020 Klagenfurt am W{\"o}rthersee, Austria
        {\tt\small  viktoriia.grushkovska@aau.at}
        \newline
$^{**}$Max Planck Institute for Dynamics of Complex Technical Systems, Magdeburg, Germany
{\tt\small zuyev@mpi-magdeburg.mpg.de}\newline
} \;and Alexander Zuyev$^{**}$}
\date{}

\begin{document}

\maketitle
\thispagestyle{empty}

\begin{abstract}
This paper focuses on the problem of constructing time-varying feedback laws that asymptotically stabilize a given part of the state variables for nonlinear control-affine systems.
It is assumed that the class of systems under consideration satisfies nonlinear controllability conditions with respect to the stabilizable variables.
Under these assumptions, a time-periodic feedback control is constructed explicitly by using the inversion of the matrix composed of the control vector fields and their Lie brackets.
The proposed control design scheme is applied to solving the leader-following problem for nonlinear multi-agent systems.
These results are illustrated with two examples of nonholonomic control problems: the partial stabilization of a rolling disc and the leader-following task for two unicycles.
\end{abstract}

\section{Introduction}
The development of control algorithms for nonholonomic mechanical systems attracts attention of many researchers because of numerous engineering applications and significant mathematical challenges.
A principal challenge in the stabilization of nonholonomic systems is caused by their underactuated structure and essentially nonlinear dynamical properties.
However, it should be noted that the requirement of asymptotic stability is redundant in many applied problems,
while only the stability with respect to some part of the state variables characterizes the desired behavior of the closed-loop system.
For example, such partial stability problems arise from the single axis stabilization of satellites~\cite{RO87,Z2001}, tracking tasks for robotic manipulators with redundant kinematics~\cite{Do05,Coh09,Mand15}, multi-agent coordination and synchronization tasks.

The stability property with respect to a part of variables was first rigorously introduced by A.~M.~Lyapunov and further studied, e.g., in~\cite{RO87,Hai95,Vo12,Zu03,Ja10,Had15,Zu15}. A detailed survey of partial stability results can be found in~\cite{Vo05}.
Despite significant progress in this area, there are only a few results on the partial stabilization of nonholonomic systems.
In particular, the paper~\cite{Hai95} examined the relation between partial stability and stability with respect to all variables for nonholonomic mechanical systems and proposed partial stability conditions. A partial stabilization problem for Lagrangian systems has been considered in~\cite{Shi00,Ko02}.
Finite-time partial stabilization problem for chained-form and cascade systems has been studied, e.g.,  in~\cite{Ja14,Chen15}. A more general class of nonlinear control system has been analyzed in~\cite{Ja06,Ja10}. The proposed sufficient partial stability conditions rely on the assumption that the system admits certain Lyapunov-like function, whose time derivative is negative definite with respect to a given part of variables. However, the problem of partial stabilization remains open for general nonholonomic systems, for which the required Lyapunov-like function cannot be effectively constructed.

In~\cite{GZ19_ES}, we have proposed practical partial asymptotic stability conditions for control-affine systems, whose averaged system has partially asymptotically stable equilibrium.
The present paper addresses the problem of explicit construction of partially stabilizing feedback laws  for nonlinear control-affine systems, whose vector fields satisfy certain Lie algebra rank condition.
We will propose a family of time-periodic feedback laws for partial stabilization of such systems and apply the obtained results to nonlinear multi-agent systems.
 Note that, although there exist  numerous results on multi-agent coordination (see, e.g.,~\cite{Ol07,Qin16} for a survey), the development of feedback control algorithms for nonholonomic agents has been addressed so far for specific systems only, such as kinematic or dynamic unicycle models~\cite{Dimos05,Aj13,Cao12,Do08,Eg01,Loi08,Yo08,Mas07,Mos07,Zhu12}, chained-form systems~\cite{Dong08}, and manipulators~\cite{Cheng08,Liu11,Tan03}.

The contribution of this paper is twofold. First, we introduce a novel approach for partial stabilization of nonlinear controllable systems by means of time-varying feedback laws.
Unlike other results on partial stabilization of nonholonomic systems, our controllers ensure exponential convergence of the trajectories (see Section 2.B).
Second, we adapt these control strategies for the leader-following problem in Section~2.C.
The stability proof is presented in Section~3, and the efficiency of the proposed controllers is illustrated by numerical simulations in Section~4.


{\it Basic notations.} Throughout this paper, we will use the following notations:

$0_n$, $0_{m\times n}$ -- $n$-dimensional zero column vector and $m\times n$-dimensional zero matrix, respectively;

$1_n$ -- $n$-dimensional column vector with all entries being equal to 1;

$I_{m\times n}$ -- $m\times n$-dimensional matrix with entries $I_{ii}=1$ and $I_{ij}=0$ whenever $i\ne j$ ($i=\overline{1,m}$, $j=\overline{1,n}$);

$\delta_{ij}$ -- Kronecker delta: $\delta_{ii}=1$ and $\delta_{ij}=0$ whenever $i\ne j$;

$\|\cdot\|$ -- the Euclidian norm in $\mathbb R^n$;

${\rm dist}(x,S)=\inf_{\xi \in S}\|x-\xi\|$ -- the distance between a point $x\in\mathbb R^{n}$ and a set $S\subset\mathbb R^{n}$;

 $B_\delta(x^*)$ -- $\delta$-neighborhood of an $x^*\in \mathbb R^n$ with $\delta>0$;

 $\partial S$, $\overline S$ -- the boundary and the closure of a set $S\subset\mathbb R^n$, respectively; $\overline S= S\cup \partial S$;

 $|S|$ -- the cardinality of a set $S$;


$ L_gf(x)$ -- the directional derivative of a vector field $f:\mathbb R^n\to\mathbb R^n $ in the direction of $g:\mathbb R^n\to\mathbb R^n $ at a point $x\in\mathbb R^n$,
 $ L_gf(x)=\lim\limits_{s\to0}\dfrac{f(x+sg(x))-f(x)}{s}$;

 $[f,g](x)$ -- the Lie bracket of vector fields $f,g:\mathbb R^n\to\mathbb R^n $ evaluated at a point $x\in\mathbb R^n$,   $[f,g](x)=L_fg(x)- L_gf(x)$.

\section{Main Results}
\subsection{Problem statement}

 Consider a nonlinear control-affine system of the form
\begin{equation}\label{sys_nonh}
  \begin{aligned}
\dot x=f_0(t,x)+\sum_{k=1}^m f_k(x) u_k,
  \end{aligned}
\end{equation}
where $x=(x_1,\dots,x_n)^\top\in D\subseteq\mathbb R^n$ is the state vector and $u=(u_1,\dots,u_m)^\top$ is the control.
 We will split the components of the state vector as $x=(y^\top,z^\top)^\top\in\mathbb R^n$ with $y\in\mathbb R^{n_1}$ and $z\in\mathbb R^{n_2}$, $n_1+n_2=n$.
 With a slight abuse of notations, the column $x$ will be also identified with $x=(y,z)$.
 The main goal of this paper is to present an explicit control strategy for stabilizing system~\eqref{sys_nonh} with respect to its $y$-variables.
For this purpose, we will exploit the notion of \emph{partial asymptotic stability}~\cite{RO87,Vo12,Zu03}.
\begin{definition}
For $y^*\in\mathbb R^{n_1}$, the set $D^*=\{x=(y,z)\in\mathbb R^n: y=y^*\}\subset D$ is \emph{$y$-asymptotically stable for the system}
\begin{equation}\label{sys}
  \dot x =f(t,x),\quad x\in\mathbb R^n,\,t\ge 0,\,f:\mathbb R^+\times \mathbb R^n \to\mathbb R^n,
\end{equation}
 if it is {\em $y$-stable} and {\em $y$-attractive}, i.e.:
\\
  $-$ \emph{$y$-stability}: for every $\Delta>0$, there exists a $\delta>0$ such that,
  for any $t_0\ge 0$ and $x(t_0)\in B_\delta(y^*) \times \mathbb R^{n_2}$, the corresponding solution $x(t)=(y(t),z(t))$ of~\eqref{sys} is uniquely defined for $t\ge t_0$, and
      $
      y(t)\in B_\Delta(y^*)\text{ for all }t\in[t_0,\infty)
      $;\\
  $-$ \emph{$y$-attractivity}: for some $\delta{>}0$ and for every $\Delta{>}0$, there exists a $t_1{\ge} 0$ such that, for any  $t_0{\ge }0$, $z(t_0)\in \mathbb R^{n_2}$,
      $$
      \text{if }y(t_0)\in B_\delta(y^*)\text{ then }y(t)\in B_\Delta(y^*)\text{ for all }t\in[t_0+t_1,\infty).
      $$
If the attractivity property holds for any $\delta{>}0$ then~$y^*$ is called to be \emph{globally  $y$-asymptotically stable for system~\eqref{sys}. }
\end{definition}
With the above definition, the main problem considered in this paper can be formulated as follows.
\begin{problem}
{Given  $y^*\in\mathbb R^{n_1}$, find a feedback law $v(t,y,z)$ such that the set $D^*=\{x=(y,z)\in\mathbb R^n: y=y^*\}$ is $y$-asymptotically stable for the closed-loop system~\eqref{sys_nonh} with $u=v(t,y,z)$.
}
\end{problem}

If $m\ge n_1$ and the matrix $F(x)=I_{{n_1}\times n} \cdot (f_1(x),f_2(x),...,f_m(x))$ is of full rank for all $x$,
then a natural approach for solving the above problem leads to defining stabilizing controls from the pseudo inversion of $F(x)$.
However, Problem~1 becomes much more challenging if the number of controls is smaller than the number of $y$-variables.
In this paper, we will present partial stabilization results for system~\eqref{sys_nonh} under certain bracket generating assumptions in case $m<n_1$.
Under these assumptions, we will extend the control design scheme from~\cite{ZuSIAM,ZGB16,GZ18} and formulate sufficient conditions for partial stabilizability of system~\eqref{sys_nonh}. Similarly to the above papers, we will define solutions of the corresponding closed-loop system in the sense of sampling.
\begin{definition}~\label{def_pi}
Given a time-varying feedback law $u=v(t,x)$, $v:{\mathbb R}^+\times D\to\mathbb R^m$, $\varepsilon>0$, and $x^0\in D$, a \emph{$\pi_\varepsilon$-solution of  system~\eqref{sys_nonh}}
corresponding to $x^0\in D$ and $v(t,x)$ is an absolutely continuous function  $x(t)\in D$, defined for $t\in {\mathbb R}^+$, such that  $x(0)=x^0$ and
$$
\begin{aligned}
\dot x(t)=f_0\big(t,x(t)\big)+\sum_{k=1}^m f_k\big(x(t)\big) v_k(t,x(\tau_j)), \quad t\in [\tau_j, \tau_j+\varepsilon),
\end{aligned}
$$
with  $\tau_j=j\varepsilon$ for $j=0,1,2,\,\dots$ .
\end{definition}

\subsection{Partially stabilizing control laws}
 In this subsection, we propose a control design scheme for partial stabilization of system~\eqref{sys_nonh} whose vector fields satisfy certain controllability condition. For the clarity of presentation, we rewrite system~\eqref{sys_nonh} as
\begin{equation}\label{sys_nonhA}
  \begin{aligned}
    \dot y=g_0(t,x)+\sum_{k=1}^{m}g_k(x)u_k,\;\;
\dot z=h_0(t,x)+\sum_{k=1}^{m}h_k(x)u_k,
  \end{aligned}
\end{equation}
where $g_k:\mathbb R^n\to\mathbb R^{n_1}$, $h_k:\mathbb R^n\to\mathbb R^{n_2}$, $g_0:\mathbb R^+\times \mathbb R^n\to\mathbb R^{n_1}$, $h_0:\mathbb R^+\times \mathbb R^n\to\mathbb R^{n_2}$,
so the vector  fields of system~\eqref{sys_nonh} are represented as
$$
\begin{aligned}
f_0(x)=\left(
         \begin{array}{c}
          {g_0(t,x)} \\
           h_0(t,x) \\
         \end{array}
       \right),\
f_k(x)=\left(
         \begin{array}{c}
           g_k(x) \\
          h_k(x) \\
         \end{array}
       \right),\quad k=\overline{1,m}.
\end{aligned}
$$
We assume that the controllability rank condition for these systems involves control vector fields and their Lie brackets. Namely, let  $D_1\subseteq\mathbb R^{n_1}$ and $D_2\subseteq\mathbb R^{n_2}$ be domains, and let $D=\{(y,z)\in\mathbb R^n:y\in D_1,z\in D_2\}.$
Assume that $f_k\in C^2(D;\mathbb R^n)$ and that the following rank condition holds for all $x\in D\subset\mathbb R^n$:
 \begin{equation}\label{rankA}
\begin{aligned}
{\rm span}\Big\{g_i(x),\, I_{n_1\times n}[f_{i_1},f_{i_2}]&(x)\Big\}={\mathbb R}^{n_1},
\end{aligned}
\end{equation}
where $i\in S_1\subseteq\{1,\dots,m\}$, $(i_1,i_2)\in S_2\subseteq\{1,\dots,m\}^2$, $|S_1|+|S_2|=n_1$. Equivalently, this means the invertibility of the following $n_1\times n_1$-matrix:
\begin{equation}\label{F}
\mathcal F(x)= \Big(\big(g_{i}(x)\big)_{i\in S_1}\ I_{n_1\times n}\big([f_{i_1},f_{i_2}]\big)_{(i_1,i_2)\in S_{2}}\Big)\;\text{for}\;\; x\in D.
\end{equation}
For this case, we adopt the family of controls~\cite{GZ18}:
\begin{equation}\label{contA}
u_k^\varepsilon(t,x) =\sum_{i\in S_1}\phi^k_i(t,x)+\dfrac{1}{\sqrt\varepsilon}\sum_{(i_1,i_2)\in S_{2}} \phi^k_{i_1i_2}(t,x),
\  k=\overline{1,m},
\end{equation}
where $\varepsilon>0$ is a small parameter, $\phi^k_i(t,x)=\delta_{ki}a_i(x)$, and $\phi^k_{i_1i_2}(t,x)$ are time-periodic functions,
\begin{equation*}
\begin{aligned}
\phi^k_{i_1i_2}&(t,x)=2\sqrt{\pi \kappa_{i_1i_2}|a_{i_1i_2}(x)|}\Big(\delta_{ki_1}\cos\Big(\dfrac{2\pi \kappa_{i_1i_2}t}{\varepsilon}\Big)+\delta_{ki_2}{\rm sign}(a_{i_1i_2}(x))\sin\Big(\dfrac{2\pi \kappa_{i_1i_2}t}{\varepsilon}\Big)\Big).
\end{aligned}
\end{equation*}
Here $\kappa_{i_1i_2}$ are pairwise distinct positive integers, and
the vector 
$$a(x)=\Big((a_{i}(x))_{i\in S_1}\ (a_{i_1i_2}(x))_{(i_1,i_2)\in S_{2}}\Big)^\top\in\mathbb R^{n_1}$$ of state-dependent coefficients
is defined as
\begin{equation}\label{aA}
a(x)=-\gamma\mathcal F^{-1}(x)(y-y^*),
\end{equation}
where $\gamma>0$, and $\mathcal F^{-1}(x)$ is the inverse matrix for $\mathcal F(x)$. The main idea behind this choice of coefficients is to ensure that the $y$-component of the solutions of~\eqref{sys_nonhA} behaves similarly to the trajectories of the system $\dot{\bar y}=-\gamma (\bar y-y^*)$ (for which, obviously, $y^*$ is the globally exponentially stable equilibrium).
%
%
To ensure that system~\eqref{sys_nonhA} can be stabilized in an arbitrary small neighborhood of the set $D^*$, we impose the following assumptions. \newpage
\begin{assumption}
{We suppose that:
\begin{itemize}
\item[A1] For any interval $I=[t_0,t_1)\subseteq {\mathbb R}^+$ and any solution $x(t)$ of system~\eqref{sys_nonhA} defined on $I$ with some control $u\in C(I;{\mathbb R}^m)$, the following property holds:
$$
\begin{aligned}
\inf\limits_{t\in I}{\rm dist}(y(t),\partial D_1){>}0\Rightarrow\inf\limits_{t\in I}{\rm dist}(z(t),\partial D_2){>}0.
\end{aligned}
$$
\item[A2] There is an  $\alpha>0$ s.t.
$\|\mathcal F^{-1}(x)\|\le \alpha$ for all $x\in D$.
\item[A3] The functions $g_k\in C^2(D;\mathbb R^{n_1})$,  $k=\overline{1,m}$, $g_0\in C^2(\mathbb R^+\times D;\mathbb R^{n_1})$, $h_0\in C(\mathbb R^+\times D;\mathbb R^{n_2})$; furthermore, the functions $g_0(t,x)$  and $h_0(t,x)$ are Lipschitz continuous with respect to $x$ uniformly in $t$, and, for any compact set $\widetilde D_1\subset D_1$, the functions $\|g_0(t,x)\|$, $\Big\|\dfrac{\partial g_0(t,x)}{\partial t}\Big\|$, $\Big\|L_{f_{j_1}}g_{j_2}(t,x)\Big\|$, and $\Big\|L_{f_0}L_{f_{k_1}}g_{k_2}(x)(t,x)\Big\|$ are bounded uniformly in $t$,  for all $y\in \widetilde  D_1$, $z\in D_2$, $t\ge0$, $j_1,j_2=\overline{0,m}$, $k_1,k_2=\overline{1,m}$.
\end{itemize}
To ensure the exponential convergence of the trajectories of system~\eqref{sys_nonhA}, we will additionally assume that
\begin{itemize}
\item [A4] for any compact set $\widetilde D_1\subset D_1$, there  exist $\alpha_l,L_g\ge 0$ ($l=\overline{0,4}$) such that
$$
\begin{aligned}
&\|g_0(t,x)-g_0(t,\tilde x)\|\le L_0\|y-\tilde y\|,\\
&\|g_0(t,x)\|\le\alpha_0\|y-y^*\|,\,\Big\|\dfrac{\partial g_0(t,x)}{\partial t}\Big\|\le\alpha_1\|y-y^*\|^\nu,\\
&\max\Big\{\Big\|L_{f_0}g_k(t,x)\Big\|,\|L_{f_k}g_0(t,x)\Big\|\Big\}\le \alpha_3\|y-y^*\|^{\nu-1},\\
&\Big\|L_{f_0}g_0(t,x)\Big\|\le \alpha_4\|y-y^*\|^\nu,\,\Big\|L_{f_0}L_{f_{i}}g_{j}(x)(t,x)\Big\|\le \alpha_5,
\end{aligned}
$$
\end{itemize}
with some $\nu\ge2$,
for all $y,\tilde y\in \widetilde  D_1$, $z,\tilde z\in D_2$, $t\ge0$, $i,j=\overline{1,n}$.}
\end{assumption}
\begin{remark}
 Assumption A1 is similar to the classical assumption on $z$-extend\-ability of solutions in partial stability theory (see, e.g.,~\cite{RO87}). For the case $D_2{=}\mathbb R^{n_2}$, this means that $z(t)$ cannot escape to infinity in finite time whenever $y(t)$ remains bounded. Such an assumption is usually satisfied for well-posed mathematical models without blow-up.
\end{remark}

The basic result of this paper is as follows.
\begin{theorem}\label{thm_partial1}
 {Let the vector fields of system~\eqref{sys_nonhA} satisfy the rank condition~\eqref{rankA}, and let Assumptions~A1--A3 hold. If the functions $u_k^\varepsilon(t,x)$ are defined by~\eqref{contA}--\eqref{aA} then,
for any $\delta>0$ and $\rho>0$, there exist  $\bar\varepsilon,\gamma,\zeta,\lambda>0$ and $t_1\ge 0$  such that, for any
$x^0\in B_\delta(y^*)\times D_2$ and $\varepsilon\in(0,\bar\varepsilon]$, the corresponding $\pi_\varepsilon$-solution of the closed-loop system~\eqref{sys_nonhA} with $u_k=u_k^\varepsilon$ and $x(0)=x^0$ is well-defined on $t\in [0,+\infty)$ and}
\begin{equation}
\begin{aligned}
 \|y(t)-y^*\|\le \zeta\|y^0-y^*\|e^{-\lambda t} +\rho&\text{ for } t\in[0,t_1),\\
 \text{and }\|y(t)-y^*\|\le \rho&\text{ for } t\in[t_1,+\infty).\\
\end{aligned}
\label{assertion1}
\end{equation}
{If, additionally, A4 is satisfied and $\gamma>\alpha_0$, then}
\begin{equation}
 \|y(t)-y^*\|= O(e^{-\lambda t})\text{ as } t\to\infty.
\label{assertion2}
\end{equation}
\end{theorem}
The proof of Theorem~\ref{thm_partial1} is given in Section~\ref{sec_proofs}.
\begin{remark}
The proposed approach can be extended to systems with higher degrees of nonholonomy. In particular, if the rank condition has the form
$$
{\rm span}\Big\{g_i(x),\, [g_{i_1},g_{i_2}](x),\, \big[g_{j_1},[g_{j_2},g_{j_3}]\big](x)\Big\}={\mathbb R}^{n_1}
$$
with $i\in S_1$, $(i_1,i_2)\in S_2$, $(j_1,j_2,j_3)\in S_3\subseteq\{1,\dots,m\}^3$, $|S_1|+|S_2|+|S_3|=n_1$, then we can take the controls from~\cite[Eqs.~(5)--(7)]{GZ18}.
\end{remark}

\subsection{Leader–following formation control}
 In this subsection, we will show an application of the obtained results to multi-agent systems.
Let us emphasize that we only present a simplified problem setup due to space limits.
Consider a system of $N+1$ heterogeneous nonholonomic agents:
\begin{subequations}
\begin{align}
  &\dot x^{\ell}=\sum_{j=1}^{m_\ell}f_j^{\ell}(x^{\ell})u_j^{\ell},\quad\ell=\overline{1,N},\label{agentA}\\
  &\dot x^{L}=f(t,x^{L}),\label{agentB}
\end{align}\label{agent}
\end{subequations}
where $x^{\ell}=(x_1^\ell,\dots,x_p^\ell)\in D_1\subseteq\mathbb R^p$, $f_j^{\ell}:D_1\to\mathbb R^p$, $u^\ell=(u_1^\ell,\dots,u_{m_\ell}^\ell)^\top\in\mathbb R^{m_\ell}$ represent the state vector, vector fields, and controls of the $\ell$-th agent, $\ell=1,2,\dots,N$; $x^L=(x_1^L,\dots,x_p^L)\in D_1$ and $f:\mathbb R^+\times D_1\to\mathbb R^p $ are the state vector and the vector field of the leader. In this section, we consider the leader-following control problem:
\begin{problem}
{Given $d_{\ell}\in {\mathbb R}^p$, $\ell=1,2,\dots,N$, the goal is to construct controls $u^{\ell}(t,x^{\ell},x^{L})$ for system~\eqref{agent} such that $\|x^{\ell}(t)-x^{L}(t)-d_{\ell}\|\to 0$ as $t\to\infty$.
}\end{problem}
Let each subsystem of~\eqref{agentA} satisfy the rank condition~\eqref{rankA}:
\begin{equation}\label{rank_agent}
\begin{aligned}
{\rm span}\Big\{f_i^\ell(x^\ell),\, [f_{i_1}^\ell,f_{i_2}^\ell]&(x^\ell)\Big\}={\mathbb R}^{p},\quad\ell=\overline{1,N},
\end{aligned}
\end{equation}
with $i\in S_1^\ell\subseteq\{1,\dots,m_\ell\}$, $(i_1,i_2)\in S_2^\ell\subseteq\{1,\dots,m_\ell\}^2$, $|S_1^\ell|+|S_2^\ell|=p$.
To apply results of the previous section,  we introduce the variables
$
y^\ell = x^\ell-x^L-d_{\ell}$, $\ell=\overline{1,N}$.
Then Problem~2 can be formulated as a partial stabilization problem for
\begin{equation}\label{agent_y}
\begin{aligned}
  &\dot y^{\ell}=\sum_{k=1}^{m_\ell}f_k^{\ell}(y^{\ell}+x^L+d^\ell)u_k^{\ell}-f(t,x^{L}),\quad\ell=\overline{1,N},\\
  &\dot x^{L}=f(t,x^{L}),
\end{aligned}
\end{equation}
with respect to the $y^\ell$-variables.
System~\eqref{agent} can be treated as a system of the type~\eqref{sys_nonhA} with
 $y=({y^1}^\top,\dots,{y^\ell}^\top)^\top\in\mathbb R^{n_1}$, $z=x^L\in\mathbb R^{n_2}$,  $n_1=Np$, $n_2=p$, $n=(N+1)p$, $m=\sum_{\ell=1}^{N}m_\ell$,
 \begin{equation*}\label{ug_agents}
 \scriptsize
 u_k=u_{k-m_{\ell_k-1}}^{\ell_k},\,
 g_k(x)=
\left(
     \begin{array}{c}
       0_{(\ell_k-1)p} \\
       f_k^{\ell_k}(y^{\ell_k}+x^L+d^{\ell_k}) \\
       0_{(N-\ell_k)p} \\
     \end{array}
   \right),\, f_k(x)=0,
 \end{equation*}
 $$
 g_0(t,x)=-1_{n_1} f(t,x^L),\, h_0(t,x)=f(t,x^L),
 $$
 where $\ell_k\in\{1,\dots, N\}$ and $m_{\ell_k}\in\{m_1,\dots,m_N\}$ are such that $k\in[m_{\ell_k}+1,m_{\ell_k}+m_{\ell_k+1}]$ $(m_0:=0)$.
Condition~\eqref{rank_agent} implies that
$$
{\rm span}\Big\{g_i(x),\, I_{n_1\times n}[g_{i_1},g_{i_2}](x)\Big\}={\mathbb R}^{n_1}
$$
with  $i\in S_1=S_1^1\cup\dots\cup S_1^N$, $(i_1,i_2)\in S_2=S_2^1\cup\dots\cup S_2^N$, $|S_1|+|S_2|=n_1$. Furthermore, matrix~\eqref{F} can be represented in the block-diagonal form:
$$
\mathcal F(x) =\left(
              \begin{array}{cccc}
                \mathcal F_1(x) & 0_{p\times p}  & \ldots & 0_{p\times p} \\
                0_{p\times p}   &\mathcal F_2(x) & \ldots & 0_{p\times p}\\
                \vdots          &  \vdots        & \ddots & \vdots \\
                0_{p\times p}   & 0_{p\times p}  & \ldots &  \mathcal F_N(x) \\
              \end{array}
            \right),
$$
$$
\begin{aligned}
\mathcal F_\ell(x)=\Big(\big(f_i^\ell(y^{\ell}{+}x^\ell&{+}d^\ell)\big)_{i\in S_1^\ell}\big([f_{i_1}^\ell,f_{i_2}^\ell](y^{\ell}{+}x^\ell{+}d^\ell)\big)_{(i_1,i_2)\in S_{2}^\ell}\Big),\;\ell=\overline{1,N}.
\end{aligned}
$$
Since $\dfrac{\partial g_k(x)}{\partial y^\ell}=\dfrac{\partial g_k(x)}{\partial x^L}$ and because of the special structure of the drift term, Assumption~1 can be rewritten as follows.
\begin{assumption}
We suppose that:
\begin{itemize}
\item[B1]
Any maximal solution $x^L(t)\in D_1$ of system~\eqref{agentB} with the initial data $x^L(t_0)\in D_1$ is defined for all $t\ge t_0$.
\item[B2] There exists an  $\alpha>0$ such that
$\|\mathcal F_\ell^{-1}(x)\|\le \alpha$ for all $x\in D_1$, $\ell=\overline{1,N}$.
\item [B3] The functions $\|f(t,x^L)\|$, $\Big\|L_{f}f(t,x^L)\Big\|$, and $\Big\|\dfrac{\partial f(t,x^L)}{\partial t}\Big\|$
are bounded for all $x^L\in D_1$,  uniformly in $t\ge0$.
\end{itemize}
\end{assumption}
Under these assumptions, we express controls~\eqref{contA} with respect to the original variables in the following form:
\begin{equation}\label{cont_agent}
\begin{aligned}
u_j^\ell=u_j^{\ell,\varepsilon}(t,x^\ell,x^L)=\sum_{i\in S_1^\ell}&\phi^{\ell,j}_i(x^\ell,x^L)+\dfrac{1}{\sqrt\varepsilon}\sum_{(i_1,i_2)\in S_{2}^\ell} \phi^{\ell,j}_{i_1i_2}(t,x^\ell,x^L),
\end{aligned}
\end{equation}
\begin{equation*}
\begin{aligned}
&\phi^{\ell,j}_i(x^\ell,x^L)=\delta_{ji}^{\ell}a_i(x^\ell,x^L),\\
& \phi^{\ell,j}_{i_1i_2}(t,x^\ell,x^L)=2\sqrt{\pi \kappa^\ell_{i_1i_2}|a^\ell_{i_1i_2}(x^\ell,x^L)|}\Big(\delta_{ji_1}{\rm sign}(a^\ell_{i_1i_2}(x^\ell,x^L))\cos\Big(\dfrac{2\pi \kappa^\ell_{i_1i_2}t}{\varepsilon}\Big)\\
&\qquad\qquad\qquad\qquad\qquad{+}\delta_{ji_2}\sin\Big(\dfrac{2\pi \kappa^\ell_{i_1i_2}t}{\varepsilon}\Big)\Big), \\
\end{aligned}
\end{equation*}
\begin{equation}\label{a_agent}
 \begin{aligned}
 a^\ell(x^\ell,x^L)=-\gamma_\ell\mathcal F_\ell^{-1}(x^\ell)(x^\ell-x^L-d^\ell),\quad\gamma_\ell>0,
\end{aligned}
\end{equation}
and the
positive integers $\kappa^\ell_{i_1i_2}$ are such that $\kappa^\ell_{i_1i_2}\ne \kappa^\ell_{j_1j_2}$ whenever $(i_1,i_2)\ne(j_1,j_2)$. Thus, each control depends only on the vector fields of the corresponding agent and its distance to the leader.
The following result is a direct consequence of Theorem~\ref{thm_partial1}.
\begin{theorem}\label{thm_agent1}
 {Let the vector fields of each subsystem of~\eqref{agentA} satisfy the rank condition~\eqref{rank_agent}, and let Assumption~2 hold.
If the functions $u_j^\ell$, $\ell=\overline{1,N}$, $j\in\overline{1,m_\ell}$, are defined by~\eqref{cont_agent}--\eqref{a_agent} then,
for any $\delta>0$ and $\rho>0$,  there exist $\bar\varepsilon,\zeta_\ell,\lambda_\ell>0$ and  $t_1^\ell\ge0$ such that, for any $\varepsilon\in(0,\bar\varepsilon]$, the $\pi_\varepsilon$-solutions $x^\ell(t)$ of each subsystem~\eqref{agentA} with the initial data $x^\ell(0)=x^{0,\ell}\in B_\delta(x^L(0))$ are well-defined on $t\in [0,+\infty)$ and}
$$
\begin{aligned}
 & \|x^\ell(t)-x^L(t)-d^\ell\|\le \zeta_\ell\|x^{0,\ell}-x^L(0)-d^\ell\|e^{-\lambda_\ell t} +\rho_\ell \text{ for } t<t_1^\ell,\\
& \|x^\ell(t)-x^L(t)-d^\ell\|\le \rho \text{ for } t\in[t_1^\ell,\infty),
\end{aligned}
 $$
 {provided that $\gamma_\ell>\dfrac{1}{\rho}\sup\limits_{t\ge0,x^L\in D_2}\|f^L(t,x^L)\|$}.
\end{theorem}

\section{Proof of Theorem~\ref{thm_partial1}}\label{sec_proofs}

The proof is conceptually similar to the proof of~\cite[Theorem~1]{GZ18}.  Below we describe its main idea and highlight the main differences caused by possibly unbounded behavior of the $z$-component of trajectories of the closed-loop system.
Throughout the proof, we define the solutions of system~\eqref{sys_nonhA} in the sense of Definition~\ref{def_pi}.
We will prove assertion~\eqref{assertion2} of Theorem~\ref{thm_partial1} and briefly outline the proof of~\eqref{assertion1} due to lack of space.

First we define several constants and sets which will be exploited throughout the proof.
Let $\delta>0$ be such that $B_{\delta}(y^*)\subset  D_1$, and let us take a  compact $\widetilde D_1$  such that $D_0=\overline{B_{\delta}(y^*)}\subset \widetilde D_1\subseteq D_1$,
so $d={\rm dist}(\partial D_0,\partial \widetilde  D_1)>0$. Assume $x^0=(y^0,z^0)\in \widetilde D=\{x=(y,z)\in \mathbb R^n:y\in D_0,z\in D_2\}$ and denote $U(x^0)=\max\limits_{0\le t\le \varepsilon}\sum_{i=1}^m |u_i^\varepsilon(t,x^0)|$.
Using H\"older's inequality and formula~\eqref{aA}, one can show that, for any $x^0\in \widetilde  D$,
\begin{equation}\label{U}
\begin{aligned}
U(x^0)&\le{c_u}\sqrt{\dfrac{\|y^0-y^*\|}{\varepsilon}},
\end{aligned}
\end{equation}
where $c_u=\sqrt{\alpha}\big(|S_1|(\alpha\gamma\|y^0-y^*\|)^{2/3}+2\sqrt{2\pi}\sum_{(i_1,i_2)\in S_2}\kappa_{i_1i_2}^{2/3}\big)^{3/4}$ is monotonically increasing with respect to  $\varepsilon$ and $\|y^0-y^*\|$.

Now we come to the main part of the proof, which we begin by showing that the solutions of system~\eqref{sys_nonhA} in $D$ are well-defined  on the time interval $[0,\varepsilon]$ with some $\varepsilon>0$. For this purpose we estimate the components of solutions of~\eqref{sys_nonhA} with $x^0\in \widetilde D$ by the following integral representations:
$$
\begin{aligned}
\|y(t)-y^0\|\le &t\Big(\alpha_0\|y^0-y^*\|+U(x^0)M_g\Big)+U(x^0)L_g\int_0^t\|z(s)-z^0\|ds\\
&+\Big(L_0+U(x^0)L_g\Big)\int_0^t\|y(s)-y^0\|ds,\\
\|z(t)-z^0\|\le &t(M_{h0}+U(x^0)M_h)+\Big(L_{h0}+U(x^0)L_{h}\Big)\int_0^t\|x(s)-x^0\|ds,
\end{aligned}
$$
where  $\alpha_0,L_0$ are defined from A4,
$$
\begin{aligned}
M_g=\sup\limits_{\hspace{-1em}\underset{1\le k\le m}{x\in \widetilde D}}\|g_k(x)\|,\; M_h=\sup\limits_{\hspace{-1em}\underset{1\le k\le m}{x\in \widetilde D}}\|h_k(x)\|,\; M_{h0}=\sup\limits_{\hspace{-1em}\underset{t\in[0,\infty)}{x\in \widetilde D}}\|h_0(t,x)\|,
\end{aligned}
$$
and $L_g,L_{h0},L_h$ are Lipschitz constants of the functions $g_0,h_0$, and $h$ in $\widetilde D$, respectively.
 After applying the Gr\"{o}nwall--Bellman inequality, the above estimates take the form:

{
\begin{align}
\|y(t)-y^0\|\le &\Big(t(\alpha_0\|y^0-y^*\|+U(x^0)M_g)\nonumber\\
&+U(x^0)L_g\int_0^t\|z(s)-z^0\|ds\Big)e^{\big(L_0+L_gU(x^0)\big)t},\label{xx0A}\\
\|z(t)-z^0\|\le &\big(t(M_{h0}+U(x^0)M_h)\nonumber\\
&+\Big(L_{h0}+U(x^0)L_{h}\Big)\int_0^t\|y(s)-y^0\|ds\big)e^{\big(L_{h0}+U(x^0)L_{h}\big)t}.\label{xx0B}
\end{align}
}
By substituting~\eqref{xx0B} into~\eqref{xx0A}, integrating by parts the appearing expressions, and applying again the Gr\"{o}nwall–Bellman inequality, we conclude that
$$
\begin{aligned}
\|y(t)-y^0\|\le t &\left(U(x^0)\left(M_g+L_gt\left(M_{h0}+U(x^0)M_h\right) e^{\left(L_{h0}+U(x^0)L_{h}\right)t}\right)+\alpha_0\|y^0-y^*\|\right)\\
&\times e^{t\left(L_0+L_gU(x^0)\left(1+\dfrac{t}{2}\left(L_{h0}+U(x^0)L_{h}\right) e^{t\left(L_0 +L_gU(x^0)+\left(L_{h0}+U(x^0)L_{h}\right)\right)}\right)\right)}.
\end{aligned}
$$
Thus, for any $\varepsilon>0$ and any $t\in[0,\varepsilon]$,
\begin{equation}\label{yy0}
\begin{aligned}
\|y(t)-y^0\|\le c_y\sqrt{\varepsilon\|y^0-y^*\|},
\end{aligned}
\end{equation}
where {
$$c_{y}=\left(c_u\left(M_g+L_g\sqrt\varepsilon\left(M_{h0}\sqrt\varepsilon+c_uM_h\right)e^{\left(L_{h0}\sqrt\varepsilon+c_uL_{h}\right)\sqrt\varepsilon}\right) +\alpha_0\sqrt{\|y^0-y^*\|}\right)e^{c_{y1}\sqrt\varepsilon},$$ $$
\begin{aligned}
c_{y1}{=}&\sqrt\varepsilon L_0{+}L_gc_u\sqrt{\|y^0{-}y^*\|} \left(1{+} \dfrac{\sqrt\varepsilon}{2}\left(L_{h0}\sqrt\varepsilon{+}c_uL_{h}\right) e^{\varepsilon L_0{+}\sqrt\varepsilon\left(L_{h0}\sqrt\varepsilon{+}c_uL_{h}\right) {+}L_gc_u\sqrt{\varepsilon\|y^0{-}y^*\|}}\right).
\end{aligned}
$$
}
Let us underline that the above parameters are monotonically increasing with respect to  $\varepsilon$ and $\|y^0-y^*\|$.
From~\eqref{yy0}, one can explicitly determine such an $\varepsilon_1>0$ that,
 for any $\varepsilon\in[0,\varepsilon_1]$, the $y$-component of the solution of system~\eqref{sys_nonhA} with controls~\eqref{contA} and initial conditions $x(0)\in\widetilde D$ remains in $\widetilde D_1$ for $t\in[0,\varepsilon_1]$. Together with the assumption A1 this implies the well-posedness of the whole solution $x(t)$  in $\widetilde D$ for $t\in[0,\varepsilon_1]$.

The next step is to analyze the value $y(\varepsilon)$. For this purpose we expand the solution of system~\eqref{sys_nonh} into the Chen--Fliess series with taking into account~\eqref{contA} and consider its $y$-component:
\begin{equation}\label{chenA}
\begin{aligned}
y(\varepsilon)=y^0+\varepsilon\mathcal F(x^0)a(x^0)+\varepsilon g_0(0,x^0)+\sigma_1(\varepsilon,x^0)&+r_0(\varepsilon)+r_1(\varepsilon),
\end{aligned}
\end{equation}
$$
\begin{aligned}
\sigma_1(\varepsilon,x^0)=&-\varepsilon\mathcal F(x^0)a(x^0)+\sum_{k=1}^mg_k(x^0)\int\limits_0^\varepsilon  u_k(s_1)d_{s_1}\\
&+\sum_{k_1,k_2=1}^{m}I_{n_1\times n}L_{f_{k_2}}f_{k_1}(x^0)\int\limits_{0}^{\varepsilon}\int\limits_0^{s_1}u_{k_1}(s_1)u_{k_2}(s_2)ds_2ds_1,\\
r_0(\varepsilon)=&\int\limits_0^\varepsilon\int\limits_0^{s_1}\Big(L_{f_0}g_0(s_2,x(s_2))+\dfrac{\partial g_0(s_2,x(s_2))}{\partial s_2}\\
&+\sum_{k=1}^m \Big( L_{f_k}g_0(s_2,x(s_2))u_k(s_2)+L_{f_0} g_k(s_2,x(s_2))u_k(s_1)\Big)ds_2ds_1,\\
r_1(\varepsilon)=&
\sum_{k_1,k_2=1}^m\int\limits_0^\varepsilon\int\limits_0^{s_1}\int\limits_0^{s_2}\Big(L_{f_0}L_{f_{k_2}}g_{k_1}(x(s_3))+\sum_{k_3=1}^{m}L_{f_{k_3}}L_{f_{k_2}}g_{k_1}(x(s_3))u_{k_3}(s_3)\Big) \\ 
& \times u_{k_1}(s_1)u_{k_2}(s_2)ds_3ds_2ds_1.
\end{aligned}
$$
Using the proposed formula for the control functions and exploiting estimate~\eqref{yy0} together with assumptions A2--A3, one can  show that there exist $c_{\sigma_1},c_0,c_1\ge0$ such that
$$
\begin{aligned}
&\|\sigma_1(\varepsilon,x^0)\|\le c_{\sigma_1}\big(\varepsilon\|y^0-y^*\|\big)^{3/2},\\
&\|r_0(\varepsilon)\|\le c_0\varepsilon^{3/2}\|y^0-y^*\|^{\nu/2},\\
&\|r_1(\varepsilon)\|\le c_1\big(\varepsilon\|y^0-y^*\|\big)^{3/2}.
\end{aligned}
$$
Assume that $\varepsilon<\dfrac{1}{\gamma}$ and $\gamma>\alpha_0$.
Exploiting the above estimates together with formula~\eqref{aA}, we deduce that
$$
\begin{aligned}
\|y(\varepsilon)-y^*\|\le &(1-\varepsilon(\gamma-\alpha_0))\|y^0-y^*\|+(c_{\sigma_1}+c_1)\big(\varepsilon\|y^0-y^*\|\big)^{3/2}+c_0\varepsilon^{3/2}\|y-y^*\|^{\nu/2}.\\
\end{aligned}
$$
For any $\lambda\in(0,\gamma-\alpha_0)$, let 
$$\varepsilon_2=\dfrac{1}{\delta^{\nu-2}}\left(\dfrac{\gamma-\alpha_0-\lambda}{(c_{\sigma_1}+c_1)\sqrt{\delta^{3-\nu}}+c_0}\right)^2.$$
 Then, for any $\varepsilon\in(0,\hat\varepsilon=\min\{\varepsilon_1,\varepsilon_2, \gamma^{-1}\})$, we have
$$
\|y(\varepsilon)-y^*\|\le(1-\varepsilon\lambda)\|y^0-y^*\|.
$$
Thus, $y(\varepsilon)\in D_0$ and, furthermore, for any $t\in[0,\varepsilon]$:
$$
\|y(t)-y^*\| \le \|y^0-y^*\|+\|y(t)-y^0\| \le \|y^0-y^*\|+ c_y\sqrt{\varepsilon\|y^0-y^*\|},
$$
what follows from the triangle inequality and estimate~\eqref{yy0}.
Summarizing all the above, we end up with the following intermediate result:  for any $\varepsilon\in(0,\hat\varepsilon)$,  the  solutions $x(t)$ of system~\eqref{sys_nonhA} with controls~\eqref{contA} and initial conditions $y^0\in B_\delta(y^*)$, $z^0\in D_2$ are well-defined in $\widetilde D$ for $t\in[0,\varepsilon_1]$ and, furthermore,
$$
\|y(\varepsilon)-y^*\|\le \|y^0-y^*\|\le\delta.
$$
Since all the parameters defined throughout the proof do not increase with decreasing $\|y^0-y^*\|$ and $\varepsilon$, the latter inequality yields that the  solutions $x(t)$ of system~\eqref{sys_nonhA} are  well-defined in $\widetilde D$ also for $t\in[0,2\varepsilon]$, and
$$
\begin{aligned}
\|y(2\varepsilon)-y^*\|\le (1-\lambda\varepsilon)^2\|y^0-y^*\|\le e^{-2\lambda\varepsilon}\|y^0-y^*\|.
\end{aligned}
$$
Iteration of  the above argumentation for $x^0\in\widetilde D$ gives the following time-periodic decay rate estimate:
$$
\|y(N\varepsilon)-y^*\|\le \|y^0-y^*\|e^{-N\lambda \varepsilon}\text{ for }N\in\mathbb N.
$$
Similarly to the proof of~\cite[Theorem~1]{GZ18}, one can show that
$$
\|y(t)-y^*\| \le \|y^0-y^*\|e^{-\lambda(t-\varepsilon)}+c_y\sqrt{\varepsilon\|y^0-y^*\|e^{-\lambda (t-\varepsilon)}}.
$$
This estimate proves assertion~\eqref{assertion2}.

Similarly, to prove assertion~\eqref{assertion1} we show that there exist $\zeta_1,\zeta_2\ge0$, $\tilde\lambda,\tilde\varepsilon>0$ such that, for all $\varepsilon\in(0,\tilde\varepsilon)$,
$$
\begin{aligned}
&\|y(\varepsilon)-y^*\|\le(1-\varepsilon\lambda)\|y^0-y^*\|+\zeta_1\varepsilon,\\
&\|y(t)-y^*\|\le \|y^0-y^*\|+\zeta_2\sqrt{\varepsilon}\text{ for }t\in[0,\varepsilon].
\end{aligned}
$$
These estimates imply that, for any $\delta,\rho>0$, there exists an $\hat\varepsilon>0$ such that, for all $\varepsilon\in(0,\hat\varepsilon)$, the $y$-component of the solutions of system~\eqref{sys_nonhA} with controls~\eqref{contA} and initial conditions $y^0\in B_\delta(y^*)$, $z^0\in D_2$, enters the $\rho$-neighborhood of $y^*$ after a finite time, and remains there.

\section{Examples}

\subsection{Position stabilization of a rolling disk}
\begin{figure*}[tpt]
\begin{center}
 \begin{minipage}{0.75\linewidth}
\begin{center}
\includegraphics[width=1\linewidth]{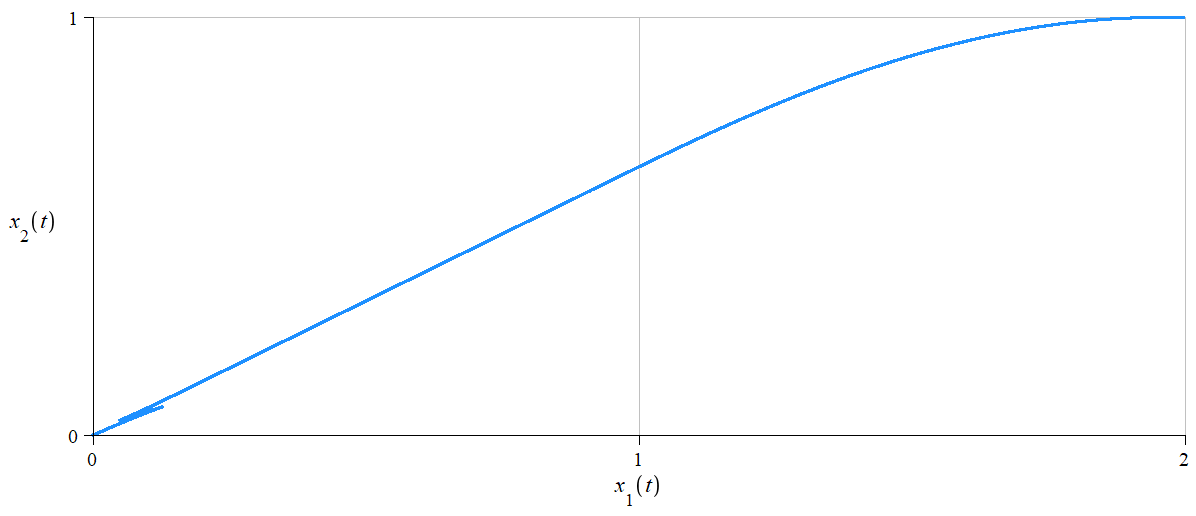}
\end{center}
 \end{minipage}\hfill
 \begin{minipage}{0.75\linewidth}
\begin{center}
\includegraphics[width=1\linewidth]{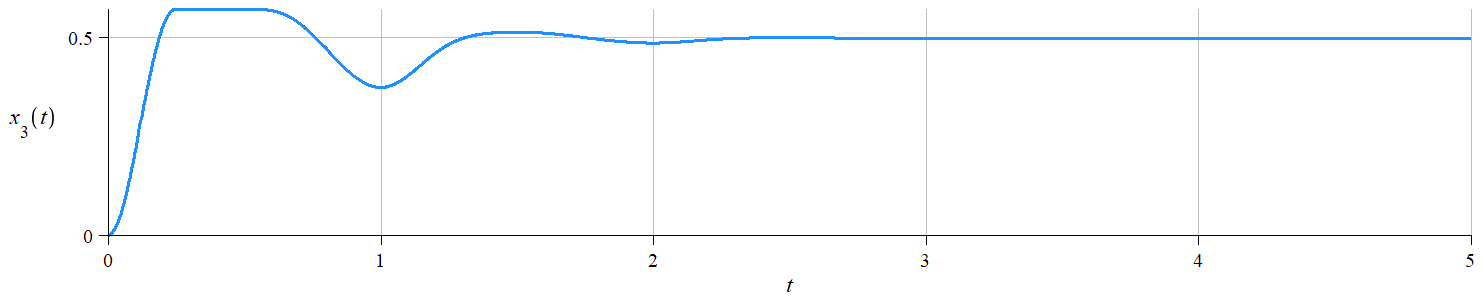}
\includegraphics[width=1\linewidth]{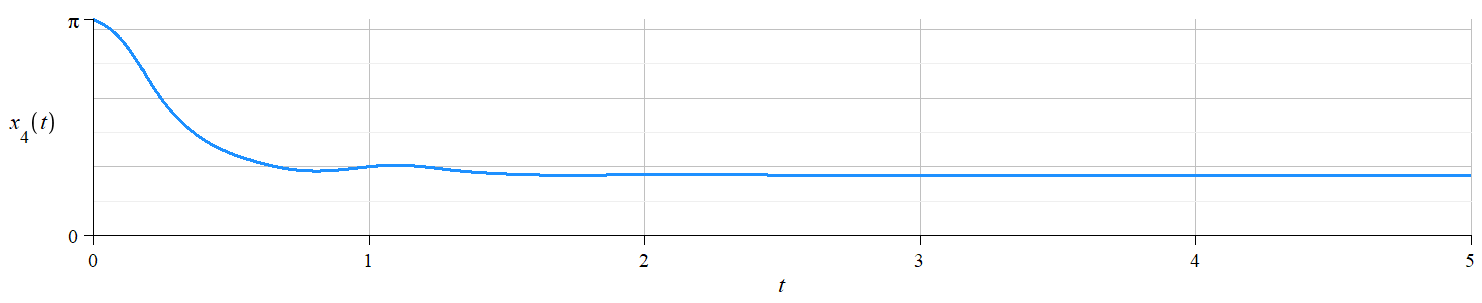}
\end{center}
 \end{minipage}
\end{center}
\caption{Projection of the trajectory of system~\eqref{disc},~\eqref{cont_disc} on the $(x_1,x_2)$-plane (left) and time-plots of $x_3(t)$ and $x_4(t)$ (right).}\label{fig_disc}
\end{figure*}

 Consider the control system describing the motion of a unit disc rolling on the horizontal plane~(cf. \cite{LC90}):
\begin{equation}\label{disc}
\dot x = f_1(x)u_1+f_2(x) u_2,\;\;x\in\mathbb R^4,\,u=(u_1,u_2)^\top\in\mathbb R^2,
\end{equation}
where $f_1(x)=\big(\cos x_3,\sin x_3,0,1\big)^\top$, $f_2(x)=\big(0,0,1,0\big)^\top$. Here $(x_1,x_2)$ are the position coordinates, the angle $x_3$  characterizes the  orientation of the disc with respect to the $x_1$-axis, and $x_4$ denotes the angle  between a fixed radial axis on the disk and the vertical axis.
The controllability condition for this system involves first- and second order Lie brackets of the vector fields $f_1$ and $f_2$:
$${\rm span}\big\{f_1(x),\,f_2(x),\,[f_1,f_2](x),\,\big[[f_1,f_2],f_2\big](x)\big\}=\mathbb R^4$$
for all $x\in\mathbb R^4$.
However, if only the position stabilization is required, then one can design stabilizing controls based on the rank condition with first order Lie brackets only. Namely, denoting $y=(x_1,x_2)^\top$, $z=(x_3,x_4)^\top$, we obtain the system of the type~\eqref{sys_nonhA} with $g_0(t,x)=h_0(t,x)\equiv 0$, $g_1(x)=\big(\cos x_3,\sin x_3\big)^\top$, $g_2(x)=\big(0,0\big)^\top$,  $h_1(x)=\big(0,1\big)^\top$, $h_2(x)=\big(1,0\big)^\top$, $n_1=n_2=2$. The rank condition~\eqref{rankA} holds with $S_1=\{1\}$, $S_{2}=\{(1,2)\}$, i.e.
$$
{\rm span}\Big\{g_1(x),\, I_{2\times 4}[f_{1},f_{2}](x)\Big\}=\mathbb R^{2}\text{ for all }x\in\mathbb R^4.
$$
To stabilize the set $D^*=\{x=(y,z)\in\mathbb R^4: y=0\}$ with respect to the $y$-variable, we apply control laws~\eqref{contA} with $\kappa_{12}=1$:
\begin{equation}\label{cont_disc}
\begin{aligned}
&u_1=a_1(x)+2\sqrt{\dfrac{\pi|a_{12}(x)|}{\varepsilon}}\cos\Big(\dfrac{2\pi\kappa_{12}t}{\varepsilon}\Big),\\
&u_2=2\sqrt{\dfrac{\pi|a_{12}(x)|}{\varepsilon}}{\rm sign}(a_{12}(x))\sin\Big(\dfrac{2\pi\kappa_{12}t}{\varepsilon}\Big),
\end{aligned}
\end{equation}
where
$a_1(x)=-\gamma\big(x_1\cos x_3+x_2\sin x_3\big)$ and $a_{12}(x)=-\gamma\big(x_1\sin x_3-x_2\cos x_3\big)$.
Figure~\ref{fig_disc} illustrates the behavior of system~\eqref{disc}--\eqref{cont_disc} with the initial condition $x(0)=(2,1,0,\pi)^\top$
and
parameters $\varepsilon=1$, $\gamma=5$. Note that, since the constructed controls tend to 0 as $x(t)\to D^*$, the solution components $x_3(t),x_4(t)$ tend to constant values.

\subsection{Leader-following strategies.}
\begin{figure*}[ht]
\begin{center}
 \begin{minipage}{0.75\linewidth}
\begin{center}
\includegraphics[width=1\linewidth]{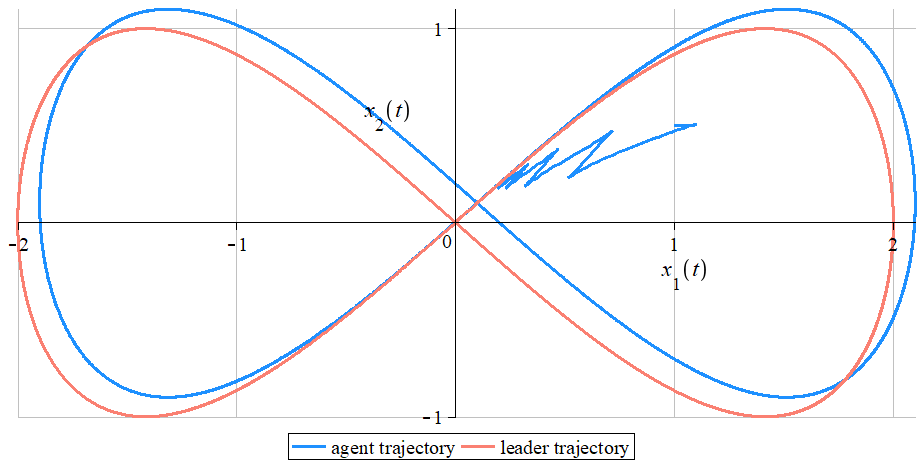}
\end{center}
 \end{minipage}\hfill
 \begin{minipage}{1\linewidth}
\begin{center}
\includegraphics[width=0.49\linewidth]{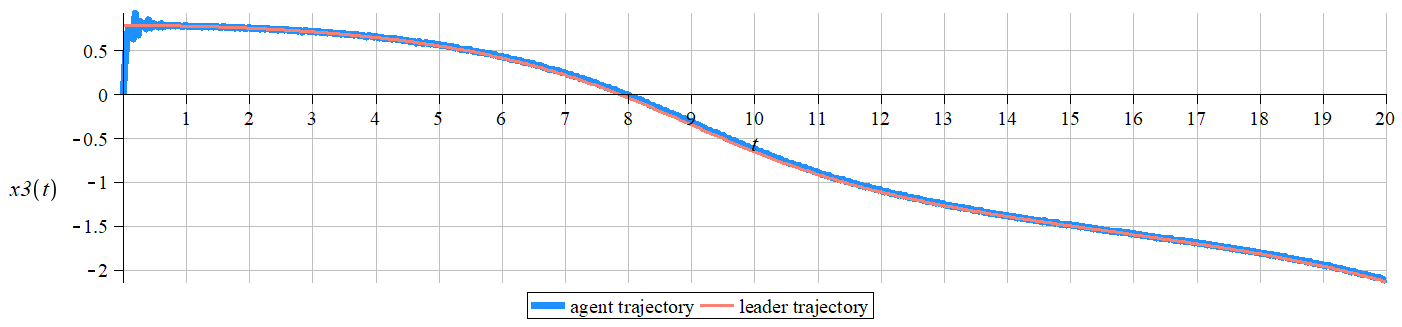}
\includegraphics[width=0.49\linewidth]{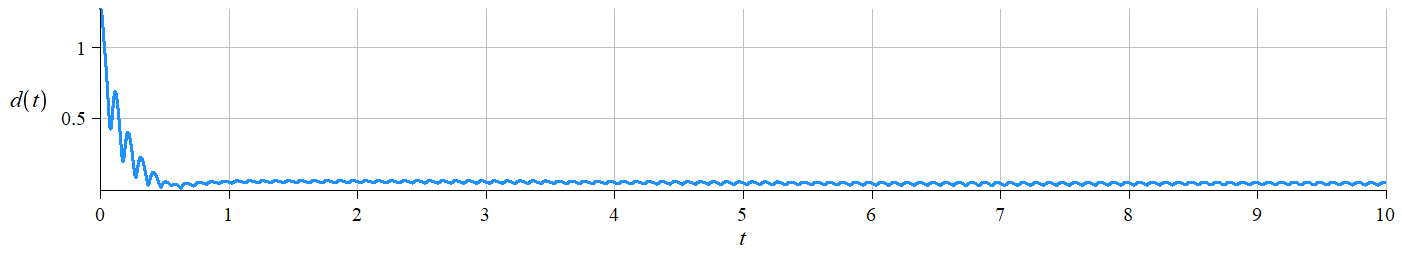}
\end{center}
 \end{minipage}
\end{center}
\caption{Projections of the trajectory of system~\eqref{uni} on the $(x_1^1,x_2^1)$- and $(x_1^L,x_2^L)$-planes (left) and time-plots of $x_3(t)$ and $d(t)=\|x^1-x^L-d^1\|$ (right).}\label{fig_uni}
\end{figure*}

To illustrate Theorem~\ref{thm_agent1}, consider Problem~2 for a system consisting of two unicycles, one of which (the leader) moves along the figure eight. The system dynamics is described by
\begin{equation}\label{uni}
\begin{aligned}
&\dot x^1 = f_1(x^1)u_1^1+f_2(x^1) u_2^1,\;\;x^1\in\mathbb R^3,\,u^1=(u_1^1,u_2^1)^\top\in\mathbb R^2,\\
&\dot x^L_1= 0.2\cos (0.1 t),\,\dot x^L_2= -0.2,\\
&\dot x^L_3=-\dfrac{0.2\sin(0.1t)\big(\cos^2(0.1t)+0.5\big)}{4\cos^4(0.1t)-3\cos^2(0.1t)+1},
\end{aligned}
\end{equation}
where
$f_1(x^1)=\big(\cos (x_3^1),\sin (x_3^1), 0\big)^\top$, $f_2(x^1)=\big(0,0,1\big)^\top$. It is easy to see that condition~\eqref{rank_agent} is satisfied with $S_1^1=\{1,2\}$, $S_2^1=\{(1,2)\}$:
$$
{\rm span}\big\{f_1(x^1),\,f_2(x^1),\,[f_1,f_2](x^1)\big\}=\mathbb R^3\text{ for all }x\in\mathbb R^3.
$$
Then we apply Theorem~\ref{thm_agent1} with
{
\begin{align}
   & u_1^1=a_1^1(x^1,x^L)+\sqrt{\dfrac{4\pi|a_{12}^1(x^1,x^L)|}{\varepsilon}}\cos\dfrac{2\pi t}{\varepsilon},\label{cont_uni}\\
   & u_2^1=a_2^1(x^1,x^L)+{\rm sign}(a_{12}^1(x^1,x^L))\sqrt{\dfrac{4\pi|a_{12}(x,\gamma)|}{\varepsilon}}\sin\dfrac{2\pi t}{\varepsilon},\nonumber
 \end{align}
 }
$$
    \begin{aligned}
  &  a_1^1(x^1,x^L)=  (x_1^1-x_1^L)\cos x_3+(x_2^1-x_2^L)\sin x_3^1,\\
 &a_2(x^1,x^L) =  x_3^1-x_3^L, \\
  &  a_{12}(x^1,x^L) =  (x_1^1-x_1^L)\sin x_3^1-(x_2^1-x_2^L)\cos x_3^1.
    \end{aligned}
$$
Figure~\ref{fig_uni} depicts the results of simulations for $\varepsilon=0.1$, $\gamma=10$, $d^1=(0.1,0.1,0)$, $x^1(0)=(1,0.5,0)^\top$, $x^L(0)=(0,0,\dfrac{\pi}{4})^\top$.


\section{Conclusions}
As one can see from Figure~\ref{fig_uni}, the proposed partially stabilizing controllers perform quite well for solving the coordination problem for nonholonomic agents.
Although we have considered rather simple formulation of the leader-following problem in Section~2.3 and Section~4, we anticipate that the results of this paper will form a basis for further treatment of the general coordination problem for multi-agent systems under essentially nonlinear controllability conditions.

\end{document}